\documentclass[11pt,sumlimits,nointlimits,namelimit,leqno]{amsart}
\usepackage{amsopn,amstext,amsfonts,amssymb,amscd,graphicx,tabularx,oldgerm,stmaryrd}
\usepackage{amsmath,amssymb,amsthm,latexsym,anysize,amssymb,amsthm}
\usepackage[mathscr]{eucal}
\usepackage[all]{xy}
\usepackage{mathrsfs}
\usepackage{paralist}
\usepackage[colorlinks]{hyperref}
\usepackage{hyperref,pdfsync}
\usepackage{CJKutf8}

 \makeatletter

\newtheorem{lemma}[equation]{Lemma}
\newtheorem*{lemma*}{Lemma}
\newtheorem{theorem}[equation]{Theorem}
\newtheorem*{theorem*}{Theorem}

\newtheorem*{corollary*}{Corollary}

\newtheorem*{proposition*}{Proposition}

\newtheorem*{conjecture*}{Conjecture}

\newtheorem*{observation*}{Observation}

\newtheorem*{claim*}{Claim}
\newtheorem*{hypothesisa*}{Hypothesis A}

\theoremstyle{definition}

\newtheorem*{definition*}{Definition}

\theoremstyle{example}

\newtheorem*{example*}{Example}
\newtheorem*{question*}{Question}

\theoremstyle{remark}

\newtheorem*{remark*}{Remark}

\def\proof{\medbreak\noindent{\scshape Proof.}\enspace}%
\def\pf{\proof}     \def\qed{\qedmark\medbreak}%
\def\qedmark{{\enspace\vrule height 6pt width 5pt depth 1.5pt}}%

\DeclareMathAlphabet\eusm{U}{eus}{m}{n}

\def\makeop#1{\expandafter\def\csname#1\endcsname
  {\mathop{\rm #1}\nolimits}\ignorespaces}
\def\doOpname#1{
 #1{ab} #1{ad} #1{Ad} #1{adj} #1{Alb} #1{Alt} #1{Ann} #1{arith}#1{Aut}  #1{Br}
#1{can} #1{Card} #1{card} #1{Char} #1{Cl} #1{Coker} #1{coker} #1{Con}
#1{Cond} #1{cond} #1{Corr} #1{curl}
#1{depth} #1{diag} #1{Disc} #1{disc} #1{Div} #1{div} #1{dom}
#1{End} #1{Ext}   #1{Fil} #1{Frac} #1{Frob} #1{Fr}
#1{Gal} #1{geom} #1{GL} #1{GSpin} #1{GSp} #1{Gr} #1{gcd}#1{gr}   #1{grad} #1{Hom} #1{hom}
#1{Id} #1{id} #1{Im} #1{im} #1{Ind} #1{ind} #1{Int} #1{inv} #1{Irr} #1{irr} #1{Isom}
#1{Ker} #1{ker} #1{lcm} #1{len} #1{length} #1{Lie} #1{Max} #1{max} #1{Min} #1{min} #1{Nr} #1{Nrd}
#1{O} #1{Ob} #1{ob} #1{opp} #1{ord}
#1{PGL} #1{Pic}  #1{Pin}#1{Proj} #1{proj} #1{Pr} #1{pr} #1{PSL}
#1{Rank} #1{rank} #1{RD} #1{Re} #1{red} #1{Res}
#1{sec} #1{sep} #1{Sgn} #1{sgn} #1{sh} #1{SL} #1{SO} #1{Sp} #1{sp}
#1{Span} #1{span}  #1{Spin} #1{spin} #1{Stab} #1{stab}
#1{std} #1{Std} #1{st} #1{St} #1{Supp} #1{supp} #1{SU} #1{Swan} #1{Sym}
#1{Tor} #1{tan} #1{tor} #1{tors} #1{torsion} #1{Trd} #1{Tr} #1{tr} #1{U} #1{val}
}
\doOpname\makeop

\begin{document}

\title{ On a problem of Sidon for polynomials over finite fields}

\author{Wentang Kuo}

\address{Department of Pure Mathematics \\
University of Waterloo \\
Waterloo, ON\\  N2L 3G1 \\
Canada}
\email{wtkuo@uwaterloo.ca}

\author{Shuntaro Yamagishi}
\address{Department of Pure Mathematics \\
University of Waterloo \\
Waterloo, ON\\  N2L 3G1 \\
Canada}
\email{syamagis@uwaterloo.ca}

\thanks{The research of the first author was supported by an NSERC
  discovery grant.}

\keywords{Sidon sets, Probabilistic number theory}
\subjclass[2010]{11K31, 11B83, 11T55}

\date{\today}

\maketitle


\begin{abstract}
Let \(\omega\) be a sequence of positive integers. Given a positive
integer \(n\), we define
\[
r_n(\omega)  = | \{ (a,b)\in \mathbb{N}\times \mathbb{N}\colon a,b \in \omega, a+b = n, 0 <a<b \}|.
\]
S. Sidon conjectured that there exists a sequence \(\omega\) such that
\(r_n(\omega) > 0\) for all \(n\) sufficiently large
and, for all \(\epsilon > 0\),
\[
\lim_{n \rightarrow \infty} \frac{r_n(\omega)}{n^{\epsilon}} = 0.
\]
P. Erd\H{o}s  proved this conjecture by showing the existence
of a sequence \(\omega\) of positive integers such that
\[
\log n \ll r_n(\omega) \ll \log n.
\]
In this paper, we prove an analogue of this conjecture in $\mathbb{F}_q[T]$,
where \(\mathbb{F}_q\) is a finite field of $q$ elements.
More precisely, let \(\omega\)
be a sequence in   \(\mathbb{F}_q[T]\). Given a polynomial
\(h\in\mathbb{F}_q[T]\), we define
\[
r_h(\omega) = |\{(f,g) \in \mathbb{F}_q[T]\times \mathbb{F}_q[T] : f,g\in  \omega, f+g =h, \deg f, \deg g \leq \deg h, f\ne g\}|.
\]
We show that there exists a sequence \(\omega\) of polynomials in \(\mathbb{F}_q [T]\)
such that
\[
\deg h  \ll r_h(\omega) \ll \deg h
\]
for \(\deg h\) sufficiently large.
\end{abstract}

\section{Introduction}\label{introduction}

In the course of investigations on Fourier series by
S. Sidon, several questions arose concerning the existence
and nature of certain positive integer sequences \(\omega\) for which
\(r_n(\omega)  = | \{ (a,b)\in \mathbb{N}\times \mathbb{N}\colon a,b \in \omega, a+b = n, 0 <a<b \}|\) is
bounded or, in some sense, exceptionally small, where \(|S|\)
denotes the cardinality of the set \(S\). In particular, he asked the
following question in 1932, known as the
Sidon Problem \cite{1}:
\begin{quote}
Does there exist a sequence \(\omega\) such that
\(r_n(\omega) > 0\) for all \(n\) sufficiently large
and, for all \(\epsilon > 0\),
\[
\lim_{n \rightarrow \infty} \frac{r_n(\omega)}{n^{\epsilon}} = 0 \text{ ?}
\]
\end{quote}

In 1954, P. Erd\H{o}s answered positively to the question
by proving the following~\cite{1}:

\begin{theorem*}[Erd\H{o}s]
There exists a sequence \(\omega\) such that
\begin{equation*}
\log n \ll r_n(\omega) \ll \log n
\end{equation*}
for all \(n\) sufficiently large.
\end{theorem*}
In other words, there exists a ``thin'' set \(\omega\) such that every
positive integer sufficiently large can be represented as a sum of two elements in
\(\omega\). On the other direction, Erd\H{o}s and R\'enyi proved the following theorem in \cite{ER}
that there exists a ``thick'' set $\omega$ such that $r_n(\omega)$ is bounded for all $n$.

\begin{theorem*}[Erd\H{o}s-R\'enyi] 
For any $\varepsilon > 0$, there exists a positive number $G = G(\varepsilon)$ and a sequence $\omega$, 
such that $r_n(\omega) < G$ for all $n$
and
\begin{equation*}
|\{ m \in \omega : m \leq n \} | > n^{\frac12 - \varepsilon}
\end{equation*}
for sufficiently large $n$.
\end{theorem*}
We note that the result is best possible up to the $\varepsilon $ term.
One way to see this fact is by the pigeon hole principle. Suppose we have $\omega_0 \subseteq \mathbb{N}$,
where $r_n(\omega_0) < G$ for all $n \in \mathbb{N}$. Given any  $m_1, m_2 \in \{ m \in \omega_0 : m \leq n \} $, we have
$1 <  m_1 + m_2 \leq 2n$. Therefore, by the pigeon hold principle, it follows that
$$
G > \max_{1 <  m \leq 2n} r_{m}(\omega_0)  \geq \frac{ | \{ m \in \omega_0 : m \leq n \} |^2 - | \{ m \in \omega_0 : m \leq n \} | }{2(2n-1)}.
$$
Consequently, we obtain
$$
| \{ m \in \omega_0 : m \leq n \} | \ll n^{1/2}.
$$
In this paper, we prove an analogue of these results
in the setting of \(\mathbb{F}_q [T]\). \\


Let $\omega$ be a sequence of polynomials in $\mathbb{F}_q[T]$.
For each \(h\in \mathbb{F}_q [T]\), we define
\[
r_h(\omega) = |\{ (f,g) \in \mathbb{F}_q[T]\times \mathbb{F}_q[T] : f,g\in  \omega,  h= f+g,  \deg f, \deg g \leq
\deg h, f\ne g\}|.
\]
Note $\deg f$ is the degree of $f \in \mathbb{F}_q[T] $ with the convention that $\deg 0 = - \infty$.
We prove the following results.
\begin{theorem}
\label{Theorem 0.1-1} There exists a sequence \(\omega\) of polynomials in \(\mathbb{F}_q [T]\)
such that
\[
\deg h  \ll r_h(\omega) \ll \deg h
\]
for \(\deg h\) sufficiently large.
\end{theorem}

On the other direction, we prove that there exists a ``thick'' set
with bounded value \(r_h(\omega)\).
We denote the elements of $\omega$ by $\omega = \{ f_i \}_{i \in \mathbb{N}}$, where
$\deg f_i \leq \deg f_j \ ( i<j) $.
\begin{theorem}
\label{Theorem 0.2-1} For each \(\epsilon> 0\), there exists a sequence
 \(\omega = \{ f_i \}\) of polynomials in \(\mathbb{F}_q [T]\) and a
positive integer \(K\) such that \(r_h(\omega) < K\) for all \(h\in \mathbb{F}_q [T]\) and
\(q^{\deg f_i}\ll i^{2+\epsilon}\).
\end{theorem}

For each \(h\in \mathbb{F}_q [T] \), we define
\[
t_h(\omega) = |\{ (f,g) \in \mathbb{F}_q[T]\times \mathbb{F}_q[T] : f,g \in \omega, h= f-g, \deg f, \deg g \leq \deg h\}|.
\]
We also prove the following variation of the existence of thick sets.
\begin{theorem}
\label{Theorem 0.3-1} For each \(\epsilon > 0\), there exists a sequence \(\omega = \{ f_i \}\) of polynomials in \(\mathbb{F}_q [T]\)
and a positive integer \(K'\) such that \(t_h(\omega) < K'\) for all
\(h\in \mathbb{F}_q [T] \) and \(q^{\deg f_i} \ll i^{2+\epsilon}\).
\end{theorem}

We prove our theorems following the methods of
Chapter III of ~\cite{2}, which utilizes the language of probability.
Roughly speaking, we set up a probability space to study
the probability of the events \(\{\omega|r_h(\omega)=d\}\) for all
non-negative integer \(d\).
Using the Borel-Cantelli lemma, we show that
the sequences satisfy the desired properties with probability $1$.
We also remark that Theorems \ref{Theorem 0.2-1} and \ref{Theorem 0.3-1} have been generalized
to $m$-fold sums and differences by K. E. Hare and the second author in \cite{HY}.\\

The organization of this paper is as follows. In Section
\ref{preliminaries},  we first review
the basic probability theory and state the  Borel-Cantelli lemma.
Next, in Section \ref{probability}, we state the equivalent
statements of our theorems and set up the probability space
used in our proof. In Section \ref{lemmas}, we establish several technical lemmas.
Finally, the remaining sections are devoted to the
proof of our main results.

\section{Preliminaries}\label{preliminaries}

We start with probability theory.  Let \(\{ X_j \}\) be a sequence of spaces
and write
\[
X = \prod_{j=0}^{\infty} X_j.
\]
Let \( \mathcal{M}_j\) be a \(\sigma\)-algebra of subsets of \(X_j\).
A measurable rectangle with respect to the sequence \(\{ \mathcal{M}_j \}\)
is defined to be a subset \(W\) of \(X\) which is representable in the
form
\[
W = \prod_{j=0}^{\infty} W_j,
\]
where \(W_j \in \mathcal{M}_j\) and \(W_j = X_j\)
except for finitely many \(j\). 
The following two theorems are standard results in probability theory,
see for example \cite[p. 123, Thm. 5]{2} and \cite[p. 135]{2} for reference.
\begin{theorem}
\label{theorem 1}\cite[p. 123, Theorem 5]{2}
Let \(\{ (X_j, \mathcal{M}_j, P_j) \}_{j \geq 0} \) be a sequence of probability spaces,
and write
\[
X = \prod_{j=0}^{\infty} X_j.
\]
Let \( \mathcal{M}\) be the minimal \(\sigma\)-algebra of subsets of \(X\)
containing every measurable rectangle with respect to the
sequence \(\{ \mathcal{M}_j \}\).
Then there exists a unique measure \(P\) on \( \mathcal{M} \) with the property
that for every non-empty measurable rectangle \(W\),
\begin{equation}
\label{eqn Thm 1-1}
P(W) = \prod_{j=0}^{\infty} P_j(W_j),
\end{equation}
where the \(W_j\) are defined by \(W = \prod_{j=0}^{\infty} W_j\),
\(W_j \in \mathcal{M}_j \ (j \geq 0)  \). 
\end{theorem}
We remark that the product in ~(\ref{eqn Thm 1-1}) is,
in essence, a finite product by the definition of measurable rectangles with respect to the
sequence \(\{ \mathcal{M}_j \}\). Furthermore, since
$$
P(X) = \prod_{j=0}^{\infty}P_j(X_j) = 1,
$$
the $\sigma$-algebra $\mathcal{M}$ in conjunction with the measure $P$ constitutes a probability space
$(X, \mathcal{M}, P)$.

\begin{theorem}\cite[p. 135, The Borel-Cantelli Lemma]{2}
Let $(X',\mathcal{M'}, P')$ be a probability space.
Let \(\{ W_{\ell} \}\) be a sequence of measurable events.
If
\[
\sum_{\ell=1}^{\infty} P'(W_{\ell}) < \infty,
\]
then, with probability \(1\), at most finite number of the events \(W_{\ell}\)
can occur; or, equivalently,
\[
P' \left( \bigcap_{i=1}^{\infty} \bigcup_{\ell=i}^{\infty} W_{\ell} \right) = 0.
\]
\end{theorem}

\section{Probability Space \((\Omega, \mathcal{M}, P)\)}\label{probability}

We let $q = p^s$ for a prime $p$, and denote $\mathbb{F}_q$ to be the finite field of $q$ elements.
Let  \(\mathbb{F}_q [T]\) be the polynomial ring over \(\mathbb{F}_q\). Let
\(\iota\) be any bijective map from \(\mathbb{Z}\cap [0,q-1]\) to \(\mathbb{F}_q\).
We label each of the polynomials in \(\mathbb{F}_q [T]\) as
follows. Let \(\mathbb{Z}_{\geq 0}\) be the set of all non-negative
integers.  For every \(N\in \mathbb{Z}_{\geq 0}\), we define
\[
p_N := \iota (c_0) + \iota (c_1) T + ... \ + \iota (c_n)T^n,
\]
where \(N = c_0 + c_1 q + ... \ + c_n q^n\) and \(0 \leq c_i < q \ (1 \leq i \leq n) \).
It is clear that this identification gives a one-to-one correspondence of sets
between \( \mathbb{Z}_{\geq 0} \) and \(\mathbb{F}_q[T]\). 
\\

We use \(\omega\) to denote a subsequence of the sequence of all polynomials in \(\mathbb{F}_q [T]\), i.e. \(p_0, p_1, p_2, p_3, ...\)
and \(\Omega\) to denote the space of all such sequences \(\omega\). By \(f \in \omega\), we mean
\(f \in \mathbb{F}_q [T]\) appears in the sequence \(\omega\). Given
\(N\in \mathbb{Z}_{\geq 0} \) and \(\omega\in \Omega \), we define
\[
r_N(\omega) = |\{(a,b)\in \mathbb{Z}_{\geq 0} \times  \mathbb{Z}_{\geq 0} \colon  p_a,p_b \in \omega, \ p_N = p_a + p_b, \ \deg p_a, \deg p_b \leq \deg p_N, \  a < b \}|,
\]
and
\[
t_N(\omega) = |\{(a,b)\in \mathbb{Z}_{\geq 0} \times  \mathbb{Z}_{\geq 0} \colon p_a,p_b \in \omega, \ p_N = p_a - p_b, \ \deg p_a, \deg p_b \leq \deg p_N \}|.
\]

We prove the following results which our main theorems, namely Theorems \ref{Theorem 0.1-1}, \ref{Theorem 0.2-1} and \ref{Theorem 0.3-1}, are consequences of.
\begin{theorem}
\label{Theorem 0.1} There exists a sequence \(\omega\) of polynomials in \(\mathbb{F}_q [T]\)
such that
\[
\log N \ll r_N(\omega) \ll \log N
\]
for \(N\) sufficiently large.
\end{theorem}
\begin{theorem}
\label{Theorem 0.2} For each \(\epsilon> 0\), there exists a sequence \(\omega = \{ p_{b_j} \}\) of polynomials in \(\mathbb{F}_q [T]\)
and a positive integer \(K_0\) such that \(r_N(\omega) < K_0\) for all \(N \in \mathbb{Z}_{\geq 0}\) and \(b_j \ll j^{2+\epsilon}\).
\end{theorem}
\begin{theorem}
\label{Theorem 0.3} For each \(\epsilon > 0\), there exists a sequence \(\omega = \{ p_{b_j} \}\) of polynomials in \(\mathbb{F}_q [T]\)
and a positive integer \(K_0'\) such that \(t_N(\omega) < K_0'\) for all \(N \in \mathbb{Z}_{\geq 0} \) and \(b_j \ll j^{2+\epsilon}\).
\end{theorem}
Since \(\deg  p_N \leq \log_q N < \deg  p_N+1\), we can easily
derive Theorems \ref{Theorem 0.1-1}, \ref{Theorem 0.2-1},
and \ref{Theorem 0.3-1}, from Theorems \ref{Theorem 0.1}, \ref{Theorem 0.2},
and \ref{Theorem 0.3}, respectively.  \\

We now prove the existence of the following probability space. The content of this theorem is essentially \cite[p. 141, Theorem 13]{2}.
\begin{theorem} \label{theorem 13}
Let
\[
\alpha_0, \alpha_1, \alpha_2, \alpha_3, ..,
\]
be real numbers satisfying \(0 \leq \alpha_i \leq 1 \ (i \geq 0) \).
Then there exists a probability space \((\Omega, \mathcal{M}, P)\) with the following
two properties:
\begin{itemize}
\item[(i)] For every non-negative integer \(m\), the event
  \(\mathfrak{B}_m = \{ \omega\in\Omega \colon p_m \in \omega \}\) is measurable
and \(P(\mathfrak{B}_m) = \alpha_m\).
\item[(ii)] The events \(\mathfrak{B}_0, \mathfrak{B}_1, \mathfrak{B}_2, ...\) are independent.
\end{itemize}
\end{theorem}

\pf Let \(Y\) be the space of two elements, \(y_0\) and \(y_1\) say.
For each sequence \(\omega\) we associate the sequence \(\{ x_j \}\)
of elements of \(Y\), defined by
\begin{eqnarray*}
x_j =
\left\{
    \begin{array}{ll}
         y_0, &\mbox{if } p_j \not \in \omega,  \\
         y_1, &\mbox{if } p_j \in \omega, \\
    \end{array}
\right.
\end{eqnarray*}
for \(j \geq 0\).
The space \(X\) consisting of all the sequences \(x = \{ x_j \}\)
is given by
\[
X = \prod_{j=0}^{\infty} X_j,
\]
where \(X_j = Y\) for \(j \geq 0\). Let \( \mathcal{M}_j = \{ \phi, \{ y_0 \}, \{ y_1 \}, X_j \}\), the non-trivial
\(\sigma\)-algebra of \(X_j\),
and let \(P_j\) be the probability measure on
\( \mathcal{M}_j\) such that \(P_j(\{ y_1\}) = \alpha_j\). \\

We apply Theorem \ref{theorem 1} to the sequence \(\{ X_j, \mathcal{M}_j, P_j \}\) of
probability spaces. In view of the one-to-one correspondence
between the elements of \(X\) and \(\Omega\), we may denote
the resulting probability space as \((\Omega, \mathcal{M}, P)\). \\

Now, we prove \((\Omega, \mathcal{M}, P)\) satisfies the two properties $(i)$ and $(ii)$.
Clearly, we have
\[
\mathfrak{B}_m = \{ \omega\in \Omega\colon p_m \in \omega \} = \prod_{j=0}^{\infty} W_j,
\]
where \(W_j = X_j\) for all \(j\) except \(j = m\) and \(W_m = \{ y_1 \}\).
Then, $(i)$ follows, because \(\mathfrak{B}_m \in \mathcal{M} \) by the definition of \( \mathcal{M} \), and
by ~(\ref{eqn Thm 1-1}) we have
\[
P(\mathfrak{B}_m) = \prod_{j=0}^{\infty} P_j(W_j) = P_m(\{ y_1 \}) = \alpha_m.
\]
For $(ii)$, we consider any finite subset of \(\{ \mathfrak{B}_j \}\), say  \(\mathfrak{B}_{j_1}, \mathfrak{B}_{j_2}, ... \ , \mathfrak{B}_{j_{\ell}}\).
Then, clearly we have
\[
\bigcap_{i=1}^{\ell} \mathfrak{B}_{j_i} = \bigl\{ \omega\in\Omega\colon p_{j_i} \in \omega \ (1 \leq i \leq  \ell) \bigr\}
= \prod_{j=0}^{\infty} W_j,
\]
where \(W_j = X_j\) for all \(j\) except \(j = j_1, ... \ , j_{\ell}\) and \(W_{j_i} = \{ y_1 \}\)
for \(1 \leq i \leq \ell\). Thus, by ~(\ref{eqn Thm 1-1}) and $(i)$ we obtain
\[
P \left( \bigcap_{i=1}^{\ell} \mathfrak{B}_{j_i} \right) = \prod_{j=0}^{\infty} P_j(W_j) = \prod_{i=1}^{\ell} P_{j_i}(\{ y_1 \})
= \prod_{i=1}^{\ell} \alpha_{j_i} = \prod_{i=1}^{\ell} P(\mathfrak{B}_{j_i}),
\]
from which $(ii)$ follows. \qed

\section{Technical Lemmas}\label{lemmas}

In this section, we prove several technical lemmas used in our
proofs. For each \(N\in \mathbb{Z}_{\geq 0}\), let \(p_N\in
\mathbb{F}_q[T]\) be as prescribed in the previous section. Define
\[
n := n(N) = \deg p_N = \lfloor \log_q N \rfloor.
\]

Suppose \(p \neq 2\).  Since \(\mathbb{F}_q = 2 \mathbb{F}_q \), we know there
exists \(p_{N_0}\) such that \(p_N = p_{N_0} + p_{N_0}\). It is clear that
\(\deg p_{N_0} = n\); therefore, \(q^n \leq N_0 < q^{n+1}\).
Since \(\mathbb{F}_q [T]\) is closed under addition, we can uniquely pair up the rest of polynomials of degree less than
or equal to \(n\) by
\[
p_N = p_{a} + p_{\widetilde{a}},
\]
where  \(a,\widetilde{a} \in  \mathbb{Z}_{\geq 0} \), \(a < \widetilde{a}\). We collect all such pairs
\((a,\widetilde{a})\) and form
\[A_N = \{a\in\mathbb{Z}_{\geq 0}\colon p_N = p_{a} + p_{\widetilde{a}}, a < \widetilde{a}, \text{ and } \deg p_a, \deg p_{\widetilde{a}} \leq n \}, \]
and
\[\widetilde{A}_N = \{\widetilde{a}\in\mathbb{Z}_{\geq 0}\colon p_N = p_{a} + p_{\widetilde{a}}, a < \widetilde{a}, \text{ and } \deg p_a, \deg p_{\widetilde{a}} \leq n\}.
\]
We have
\(| A_N | = | \widetilde{A}_N | = (q^{n+1}-1)/2, \)
and
\[
\{0, ... \ , q^{n+1} - 1 \} = A_N \bigcup \widetilde{A}_N \bigcup  \{N_0\},
\]
where all the unions are disjoint.
Further, it is easy to see that
\(\{ 0, 1, ... \ , q^n - 1 \} \subseteq A_N\), because
if \(0 \leq a < q^n\), then
\(p_{a}\) has degree at most \(n-1\). Thus, the corresponding
\(p_{\widetilde{a}}\) must have degree \(n\); therefore,
\( q^n \leq \widetilde{a} < q^{n+1}\).
Hence, it follows that
\begin{equation}
\label{On tilde A_N}
\widetilde{A}_N \subseteq \{ q^n, q^n + 1, ... \ , q^{n+1} - 1 \}.
\end{equation}

Let \(M:= M(N) = (q^{n+1}-1)/2\). For convenience we label the \(M\) elements of
\(A_N\) by \(a_i\), where \(1 \leq i \leq M\),
and the corresponding elements of \(\widetilde{A}_N\) by \(\widetilde{a}_i\).

We also define \(\lambda_N\) and \(\lambda '_N\) to be
\[
\lambda_N = \sum_{1 \leq i \leq M } \alpha_{a_i} \alpha_{\widetilde{a}_i},
\]
and
\[
\lambda'_N = \sum_{1 \leq i \leq M } \frac{ \alpha_{a_i} \alpha_{\widetilde{a}_i} }{ 1 - \alpha_{a_i} \alpha_{\widetilde{a}_i}}.
\]
Note when \(p=2\), for $N>0$, we do not have to consider the polynomial \(p_{N_0}\) as above. Thus we let \(M:= M(N) = q^{n+1}/2\)
and we can argue in a similar manner.   \\

Define
\[
s^*_N(\omega) = \sum_{m=0}^{N} \mathbf{1}_{\mathfrak{B}_m}(\omega),
\]
where \(\mathbf{1}_{\mathfrak{B}_m}\) is the characteristic function on the set \(\mathfrak{B}_m\).
Let \(E(f)\) denote the expectation of a random variable \(f\),
defined by \(E(f) = \int_{X} f \ dP\).
We define
\[
m^*_N = E(s^*_N) = \sum_{m=0}^{N} \alpha_m.
\]
\\

We need our sequence \(\{\alpha _j\}\) to satisfy the following condition.
\begin{hypothesisa*}
The sequence \(\{ \alpha_j \}\) of probabilities (introduced in Theorem \ref{theorem 13})
satisfies the conditions: \(0 < \alpha_j < 1 \ (j \geq 0) \), \(\{ \alpha_j \}\) is monotonic and decreasing from some point onward (i.e.
for \(j \geq j_1\)), and \(\alpha_j \rightarrow 0\) as \(j \rightarrow \infty\).
\end{hypothesisa*}
We have the following result for \(s^*_N(\omega) \) and its expected value \(m^*_N\).
\begin{lemma}
\label{Lemma 10} If, in addition to Hypothesis A,
\begin{equation}
\label{eqn lemma 10-1}
m^*_N \rightarrow \infty
\end{equation}
as \(N \rightarrow \infty\), and
\begin{equation}
\label{eqn lemma 10-2}
\sum_{N=0}^{\infty} \frac{\alpha_N}{(m^*_N)^2} < \infty,
\end{equation}
then with probability \(1\), we have  
\(s^*_N (\omega) \sim m^*_N\) as \(N \rightarrow \infty\).
\end{lemma}
\pf
We denote \(D^2(f)\) to be the variance of a random variable \(f\),
defined by
\[
D^2(f) = E\bigl( ( f - E(f) )^2 \bigr).
\]
The proof is basically an application of a variant of the strong law of large numbers 
\cite[p. 140, Theorem 11]{2}, which is as follows.
Let \(\{f_j\}\) be a sequence of independent random
variables, and let
\[
s_i (\omega) = \sum_{j=0}^{i} f_j (\omega) \ \ (i \geq 0).
\]
Suppose
we have
\[
E(f_j) > 0 \ \ \ (j \geq 0),
\]
\[
\lim_{i \rightarrow \infty} E(s_i) = \infty,
\]
and
\[
\sum_{i=0}^{\infty} \frac{D^2(f_i)}{( E(s_i) )^2 } < \infty.
\]
Then, with probability \(1\), we have
\[
s_i(\omega) =\bigl(1 + \mathbf{o}(1) \bigr) E(s_i)
\]
as \(i \rightarrow \infty\).
We know that the sets $\mathfrak{B}_j$ are independent, which is equivalent to $\mathbf{1}_{\mathfrak{B}_j}(\omega)$
being independent. Thus we apply this theorem with \(f_j(\omega) = \mathbf{1}_{\mathfrak{B}_j}(\omega)\),
and obtain our result. \qed

For every \(N,d\in \mathbb{Z}_{\geq 0}\),
we define the event \(\mathfrak{e}(N,d)\) as
\[
\mathfrak{e}(N,d) =\{\omega \in\Omega  : r_N(\omega) = d \}.
\]

As mentioned in Section \ref{introduction}, we need to study
the probability of the event \(\mathfrak{e}(N,d)\). We
start with the following lemma.

\begin{lemma}
\label{Lemma 12} For all non-negative integers \(N\) and \(d\), we have
\begin{equation}
\label{eqn lemma 12-1}
P\bigl(\mathfrak{e}(N,d)\bigr)
= \left( \prod_{1 \leq k \leq M } (1 - \alpha_{a_k} \alpha_{\widetilde{a}_k}) \right) \widetilde{\sigma}_{d}(N),
\end{equation}
where \(\widetilde{\sigma}_{0}(N) = 1\) and, if \(d \geq 1\),
\begin{equation}
\widetilde{\sigma}_{d}(N) = \sum_{1 \leq k_1 < ... < k_d \leq M } \
\prod_{1 \leq i \leq d} \frac{\alpha_{a_{k_i}} \alpha_{\widetilde{a}_{k_i}} }{1 - \alpha_{a_{k_i}} \alpha_{\widetilde{a}_{k_i}} }.
\end{equation}
\end{lemma}

\pf We begin with the case \(d=0\). It is easy to see that
\[
\mathfrak{e}(N,0) =  \bigcap_{1 \leq k \leq M} \left( \mathfrak{B}_{a_k} \cap \mathfrak{B}_{\widetilde{a}_k} \right)^{\mathbf{c}},
\]
where \(\mathbf{c}\) denotes taking the complement of the set. Since the sets \(  \mathfrak{B}_j \ (j \geq 0)  \) are
independent, we know that \(\mathfrak{B}_{a_k} \cap  \mathfrak{B}_{\widetilde{a}_k} \ (1 \leq k \leq M) \) are independent as
\(\{ a_k : 1 \leq k \leq M \} \bigcap  \{ \widetilde{a}_k : 1 \leq k \leq M \} = \varnothing\).
Thus, it follows that \( \left( \mathfrak{B}_{a_k} \cap  \mathfrak{B}_{\widetilde{a}_k} \right)^{\mathbf{c}} \ (1 \leq k \leq M)  \) are also independent.
Hence, we have
\[
P(\mathfrak{e}(N,0)) =  \prod_{1 \leq k \leq M} P \left( \left( \mathfrak{B}_{a_k} \cap \mathfrak{B}_{\widetilde{a}_k} \right)^{\mathbf{c}} \right)
= \prod_{1 \leq k \leq M} (1 - \alpha_{a_k} \alpha_{\widetilde{a}_k}).
\]

Suppose \(1 \leq d \leq M\) and \(\omega' \in \mathfrak{e}(N,d)\).
Then there exist \(k_1, k_2, ... \ , k_d\) such that \(1 \leq k_i \leq M\),
\(a_{k_i}, \widetilde{a}_{k_i} \in \omega'\)  \( ( 1 \leq i \leq  d )\), and
further, if $k \not = k_i $ and
\(1 \leq k \leq M\), then we have either $a_k \not \in \omega'$ or $\widetilde{a}_k \not \in \omega'$.
From this observation, we can deduce that
\[
P\bigl(\mathfrak{e}(N,d)\bigr) = \sum_{1 \leq k_1 < ... < k_d \leq M} P\bigl(\mathfrak{E}(k_1, ... \ ,k_d )\bigr),
\]
where \(\mathfrak{E}(k_1, ... \ ,k_d )\) is the event
\[
\bigcap_{1 \leq i \leq d} \left( \mathfrak{B}_{a_{k_i}} \cap \mathfrak{B}_{\widetilde{a}_{k_i}} \right) \ \bigcap \bigcap_{ \stackrel{1 \leq k \leq M}{k \not = k_i   (1 \leq i \leq  d)  } } \left( \mathfrak{B}_{ a_{k} } \cap \mathfrak{B}_{\widetilde{a}_{k}} \right)^{\mathbf{c}}.
\]
Again, by independence, we have
\begin{eqnarray}
P\bigl(\mathfrak{E}(k_1, ... \ ,k_d )\bigr) &=& \prod_{1 \leq i \leq d} P \left( \mathfrak{B}_{a_{k_i}} \cap \mathfrak{B}_{\widetilde{a}_{k_i}} \right) \ \cdot  \prod_{ \stackrel{1 \leq k \leq M}{k \not = k_i   (1 \leq i \leq  d)  } } P\left( \left( \mathfrak{B}_{ a_{k} } \cap \mathfrak{B}_{\widetilde{a}_{k}} \right)^{\mathbf{c}} \right)
\notag
\\
&=& \prod_{1 \leq i \leq d} \alpha_{ a_{k_i} } \alpha_{ \widetilde{a}_{k_i} } \ \cdot
\prod_{ \stackrel{1 \leq k \leq M }{k \not = k_i  (1 \leq i \leq  d)  } } \left( 1 - \alpha_{ a_{k} } \alpha_{ \widetilde{a}_{k} } \right)
\notag
\\
&=& \prod_{1 \leq k \leq M } (1 - \alpha_{a_k} \alpha_{ \widetilde{a}_{k} }) \ \cdot
\prod_{ 1 \leq i \leq d }\frac{ \alpha_{ a_{k_i} } \alpha_{ \widetilde{a}_{k_i} } }{1 - \alpha_{ a_{k_i} } \alpha_{ \widetilde{a}_{k_i} }},
\notag
\end{eqnarray}
from which the desired result follows. \\

Finally, if \(d > M\), then the sum \(\widetilde{\sigma}_{d}(N)\) is empty,
and both sides of ~(\ref{eqn lemma 12-1}) are \(0\). \qed

To estimate \(\widetilde{\sigma}_{d}(N)\), we use the following result
for  elementary symmetric functions.

\begin{lemma}\cite[p. 147, Lemma 13]{2} 
\label{Lemma 13} Let \(y_1, y_2, ... \ y_{M'}\) be \(M'\) non-negative real numbers.
For each positive integer \(d\), not exceeding \(M'\), let
\[
\sigma_d = \sum_{1 \leq  k_1 < ... \ < k_d \leq M'} y_{k_1} y_{k_2} ... \ y_{k_d},
\]
so that \(\sigma_d\) is the \(d\)-th elementary symmetric function of the \(y_k\)'s.
Then, for each \(d\), we have
\begin{equation}
\label{eqn lemma 13-3}
\frac{1}{d!} \sigma_1^d \left( 1- {d \choose 2} \frac{1}{\sigma_1^2}  \sum_{k=1}^{M'} y_k^2 \right) \leq \sigma_d \leq \frac{1}{d!} \sigma_1^d,
\end{equation}
where we interpret \(d \choose 2\) to be \(0\) when \(d=1\).
\end{lemma}


The next lemma gives us bounds on the probability of the event \(\mathfrak{e}(N,d)\) in terms of \(\lambda _N\) and \(\lambda
_N'\).
\begin{lemma}
\label{Lemma 14} Let \(N\) and \(d\) be non-negative
integers. Then we have
\begin{equation}
\label{eqn lemma 14-1}
P\bigl(\mathfrak{e}(N,d)\bigr) \leq \frac{(\lambda '_N)^d}{d!}e^{-\lambda_N}.
\end{equation}
Furthermore, if \(d \leq M\), we have
\begin{equation}
\label{eqn lemma 14-2}
P\bigl(\mathfrak{e}(N,d)\bigr)
\geq
\frac{(\lambda'_N)^d}{d!}e^{-\lambda'_N} \left( 1 - {d \choose 2} (\lambda'_N)^{-2} Q^* \right),
\end{equation}
where
\[
Q^* = \sum_{1 \leq k \leq M} \left( \frac{ \alpha_{a_k} \alpha_{ \widetilde{a}_k } }{1 - \alpha_{a_k} \alpha_{ \widetilde{a}_k } } \right)^2,
\]
and \(d \choose 2\) is interpreted to be \(0\) if \(d<2\).
\end{lemma}

\pf We note that if \(d > M\), then the event \(\mathfrak{e}(N,d)\) is empty and
~(\ref{eqn lemma 14-1}) is trivial. Suppose \(1 \leq d \leq M\). We apply
~(\ref{eqn lemma 13-3}) with $M'=M$ and \(y_k = \alpha_{a_k} \alpha_{ \widetilde{a}_k }/ (1 - \alpha_{a_k} \alpha_{ \widetilde{a}_k })\)
to estimate \(\widetilde{\sigma}_{d}(N)\) in ~(\ref{eqn lemma 12-1}); thus noting that \(\widetilde{\sigma}_{1}(N) = \lambda'_N\), we obtain
\[
P\bigl(\mathfrak{e}(N,d)\bigr)
\leq \left( \prod_{1 \leq k \leq M } (1 - \alpha_{a_k} \alpha_{ \widetilde{a}_k } ) \right) \frac{(\lambda'_N)^d}{d!},
\]
and
\[
P\bigl(\mathfrak{e}(N,d)\bigr)
\geq \left( \prod_{1 \leq k \leq M } (1 - \alpha_{a_k} \alpha_{ \widetilde{a}_k }) \right) \frac{(\lambda'_N)^d}{d!} \left( 1 - {d \choose 2} (\lambda'_N)^{-2} Q^* \right).
\]

Applying the inequality \(e^{-t/(1-t)} < 1 - t < e^{-t}\), which holds for \(0 < t < 1\), with $t = \alpha_{a_k} \alpha_{ \widetilde{a}_k } \ (1 \leq k \leq M)$, we obtain
\begin{equation}
\label{basic inequality 1}
e^{-\lambda'_N} < \left( \prod_{1 \leq k \leq M } (1 - \alpha_{a_k} \alpha_{ \widetilde{a}_k }) \right) < e^{-\lambda_N},
\end{equation}
and our result follows.
When $d=0$, we have
$$
P\bigl(\mathfrak{e}(N,d)\bigr) =  \prod_{1 \leq k \leq M } (1 - \alpha_{a_k} \alpha_{ \widetilde{a}_k }),
$$
and the result is immediate from ~(\ref{basic inequality 1}) in this case.
\qed

To estimate \(\lambda_N\) and \(\lambda _N'\), we first prove the
following lemma.

\begin{lemma}
\label{Lemma 15} If Hypothesis A is satisfied, then
\begin{equation}
\lambda '_N \sim \lambda_N
\end{equation}
as \(N \rightarrow \infty\).
\end{lemma}

\pf Recall from ~(\ref{On tilde A_N}) that if \(1 \leq k \leq M\),
then  \(q^n \leq \widetilde{a}_k < q^{n+1}\). Consequently, we have
\[
\alpha_{a_k}\alpha_{\widetilde{a}_k} < \alpha_{\widetilde{a}_k} \leq \alpha_{q^n} =  \mathbf{o}(1)
\]
as \(N \rightarrow \infty\). Therefore, we obtain
\begin{eqnarray*}
\lambda '_N - \lambda_N &=& \sum_{1 \leq k \leq M}  \alpha_{a_k}\alpha_{\widetilde{a}_k} \left( \frac{1}{1-\alpha_{a_k}\alpha_{\widetilde{a}_k}  } -1 \right)
\notag
\\
&=& \sum_{1 \leq k \leq M}  \alpha_{a_k}\alpha_{\widetilde{a}_k} \left( \frac{\alpha_{a_k}\alpha_{\widetilde{a}_k}}{1-\alpha_{a_k}\alpha_{\widetilde{a}_k}} \right)
\notag
\\
&\leq &\alpha_{q^n} \lambda '_N,
\end{eqnarray*}
from which the result follows. \qed

Therefore, it is enough to estimate \(\lambda _N\). The following
lemma gives us an estimate for \(\lambda _N\) sufficient for our purpose.

\begin{lemma}\label{Lemma 11} Suppose that the sequence \(\{ \alpha_j \}\), introduced in Theorem \ref{theorem 13},
is such that
\[
\alpha_j = \alpha \frac{(\log j)^{c'}}{j^c}
\]
for \(j \geq j_0\); where \(j_0, \alpha, c, c'\) are constants such that \(\alpha > 0\),
\(0 < \alpha_j < 1 \ (j \geq 0) \), \(0 < c < 1\), and \(c' \geq 0\).
Then, for sufficiently large $H$, there exist positive constants \(D_1\) and \(D_2\), which depend at most on
\(c\), \(c'\), and \(q\), 
such that
\begin{equation}
\label{eqn lemma 11}
\alpha^2 D_1 (\log N)^{2 c'} q^{n(1-2c)} < \lambda_N < \alpha^2 D_2 (\log N)^{2 c'} q^{n(1-2c)}
\end{equation}
for all $N > H.$

Furthermore, we have
\begin{equation}
\label{eqn lemma 11 statement 2}
m^*_N \sim \frac{\alpha}{1-c} (\log N)^{c'} N^{(1-c)}.
\end{equation}

Finally, if \(c' = 0\), then with probability \(1\), the numbers \(b_j\) of
the sequence \(\omega = \{ p_{b_j} \}\) satisfy
\begin{equation}
\label{eqn lemma 11 statement 3}
b_j \sim \left( \frac{1-c}{\alpha} j \right)^{1/(1-c)}
\end{equation}
as \(j \rightarrow \infty\).
\end{lemma}

\pf We begin by finding a lower bound for \(\lambda_N\). We assume \(p\neq 2\). The case \(p=2\) can be treated in a similar manner.
Suppose
\(N > q\cdot (j_0+1)\) from which it follows that $q^n > j_0$.  Let \(C_0\) and \(C_0'\) be the positive constants defined by
\[
C_0=\sum_{1 \leq j < j_0} \alpha_{j} ,
\]
and
\[
C_0'= \sum_{1 \leq j < j_0} 
\frac{(\log j)^{c'}}{j^c}.
\]

Since \(q^n \leq \widetilde{a}_i < q^{n+1}\), \(0 \leq a_i < q^{n+1}\),
and \((\log x)^{c'}/x^c\) is a decreasing function, we obtain
\[
\begin{split}
 \lambda_N = \sum_{1 \leq i \leq M} \alpha_{a_i}
 \alpha_{\widetilde{a}_i}& > \frac{ \alpha (\log q^{n+1})^{c'}
 }{q^{(n+1)c}}\sum_{1 \leq i \leq M} \alpha_{a_i}  \\
 &   >
\frac{\alpha^2 (\log q^{n+1})^{c'} }{q^{(n+1)c}}  \left(\sum_{j =
  (q^{n+1}-1)/2 }^{q^{n+1}-1} \frac{(\log j)^{c'}}{j^c}-C_0'\right)\\
 & > \frac{\alpha^2 (\log q^{n})^{2c'} }{q^{(n+1)c}}   \left(\sum_{j =
     (q^{n+1}-1)/2 }^{q^{n+1}-1}
   \frac{1}{j^c}-C_0'\right) .
\end{split}
\]

We know that for all \(s,t\in \mathbb{N}\), \(0<s<t\),

\begin{eqnarray}
\label{eqn lemma 11 ineq}
\frac{1}{1-c} (t+1)^{1-c} - \frac{1}{1-c} s^{1-c} 
\leq
\sum_{s \leq j \leq t} \frac{1}{j^c} 
\leq \frac{1}{1-c} t^{1-c} - \frac{1}{1-c} (s-1)^{1-c}.
\end{eqnarray}

Thus, by ~(\ref{eqn lemma 11 ineq}) we can give the following lower bound for \(\lambda_N\),
\begin{eqnarray}
\label{eqn lemma 11-1}
\lambda_N
&>&
\frac{\alpha^2 (\log q^n)^{2c'} }{q^{(n+1)c}} \left(
\frac{1}{1-c} (q^{n+1})^{1-c} - \frac{1}{1-c} \left( \frac{q^{n+1}-1}{ 2} \right)^{1-c} -C_0'\right)
\\
&=& \frac{\alpha^2 (\log q^n)^{2c'} }{(1-c)q^c} q^{n(1 - 2c)} \left(
q^{1-c} - \left(\frac{q}{2} - \frac{1}{2q^n} \right)^{1-c} -\frac{1-c}{q^{n(1-c)}}\cdot C_0'\right).
\notag
\end{eqnarray}
Since \(q^n \leq N < q^{n+1}\), we have \(\log N (1 - \log q/ \log N) < \log q^n\).
It follows from ~(\ref{eqn lemma 11-1}) that by taking $H$ sufficiently large, we obtain
\begin{equation}
\label{estimate D_1}
\lambda_N > \alpha^2 \frac{q^{1-2c}}{2(1-c)} \left( 1 - \frac{1}{2^{1-c}}  \right) (\log N)^{2c'} q^{n(1 - 2c)}
\end{equation}
for all $N>H$.

Next, we would like to find an upper bound for \(\lambda_N\).
Again, since \(q^n \leq \widetilde{a}_i < q^{n+1}\) and \(0 \leq a_i < q^{n+1}\),
by similar calculations as before we have
\[
\lambda_N < \frac{\alpha ( \log q^{n+1})^{c'}}{q^{nc}} \sum_{1 \leq i \leq M } \alpha_{a_i}
<
\frac{\alpha ( \log q^{n+1})^{c'}}{q^{nc}} \left( C_0 + \alpha \sum_{j=1}^{M} \frac{ (\log q^{n+1})^{c'} }{j^c} \right).
\]

Thus, by applying ~(\ref{eqn lemma 11 ineq}), we obtain
\begin{eqnarray*}
\lambda_N
&<&
\frac{\alpha^2 ( \log q^{n+1})^{2c'}}{q^{nc}} \left( \frac{ C_0 }{ \alpha (\log q^{n+1})^{c'} }  +
\frac{1}{1-c} \left(  \frac{q^{n+1}- 1}{2}   \right)^{1-c} \right)
\notag
\\
&=& \frac{\alpha^2 ( \log q^{n+1})^{2c'}}{(1-c)} q^{n(1-2c)}
\left(
\frac{C_0(1-c)}{ \alpha (\log q^{n+1})^{c'} q^{n(1-c)} } +
\left(
\frac{q}{2} - \frac{1}{2 q^n} \right)^{1-c}
\right).
\end{eqnarray*}
Therefore, by taking $H$ sufficiently large, we obtain that
\[
\lambda_N < \alpha^2  \frac{2^{2 c' + 1}}{(1-c)} \left( \frac{q}{2} \right)^{1-c}  (\log q^n)^{2c'} q^{n(1 - 2c)}
\]
for all $N>H.$
Then, since \(q^n \leq N < q^{n+1}\), we are done with the first part
of the lemma.  \\


Clearly, we have
\begin{eqnarray*}
m^*_N
= \sum_{j=1}^{N} \alpha \frac{(\log j)^{c'}}{j^c} + \mathbf{O}(1)
= \bigl( 1 + \mathbf{o}(1)\bigr)  \frac{\alpha}{1-c} (\log N)^{c'} N^{(1-c)},
\end{eqnarray*}
and this proves ~(\ref{eqn lemma 11 statement 2}). We note ~(\ref{eqn lemma 11 statement 2}) shows that
~(\ref{eqn lemma 10-1})
and ~(\ref{eqn lemma 10-2}) are satisfied. \\

The final assertion of the lemma follows from ~(\ref{eqn lemma 11 statement 2}),
in view of Lemma \ref{Lemma 10}, and the fact that \(s^*_{b_j} (\omega) = j\) for $ \omega = \{ p_{b_j} \}_{j \in \mathbb{N}} $; for in this way it follows that,
if \(c'=0\), we have with probability \(1\),
\[
j = s^*_{b_j} (\omega) \sim m^*_{b_j} \sim \frac{\alpha}{1-c} b_j^{1-c},
\]
or equivalently, \(b_j \sim \left( \frac{1-c}{\alpha} j
\right)^{1/(1-c)}\). \qed

We also make use of the following lemma.
\begin{lemma}\cite[p. 149, Lemma 17]{2} 
\label{Lemma 17} If \(0 < \xi \leq U\), then
\begin{equation*}
\sum_{d \geq U} \frac{\xi^d}{d!} \leq \left( \frac{e \xi}{U} \right)^U,
\end{equation*}
and if \(0 < V \leq \xi\), then
\begin{equation*}
\sum_{0 \leq d \leq V} \frac{\xi^d}{d!} \leq \left( \frac{e \xi}{V} \right)^V.
\end{equation*}
\end{lemma}


\section{Proof of Theorem \ref{Theorem 0.1}}

Let \(c = c' = 1/2\). We choose a number \(\alpha > 0\) to satisfy
$$
\alpha^2 \frac{q^{1-2c}}{2(1-c)} \left( 1 - \frac{1}{2^{1-c}}  \right) > 1.
$$
We then define a sequence \(\{ \alpha_j \}\) by
\begin{equation}
\label{eqn T0}
\alpha_j = \alpha \left( \frac{\log j}{j} \right)^{1/2}
\end{equation}
for all \(j \geq j_0\), where \(j_0\) is a positive integer
sufficient large such that the expression in
~(\ref{eqn T0}) is less than \(1/2\) for all \(j \geq j_0\).
For \(1 \leq j < j_0\), we let \(\alpha_j = 1/2\).
The precise value of \(\alpha_j\) for small \(j\) is unimportant,
but the above choices ensure \(0 < \alpha_j < 1\),
so that Hypothesis A is satisfied. By ~(\ref{estimate D_1}) in the proof of Lemma \ref{Lemma 11},
we have for all \(N\) sufficiently large that
\begin{equation}
\label{eqn T1}
\lambda_N \geq \alpha^2 D_1 \log N > \log N.
\end{equation}
Hence, we know
there exists \(\delta > 0\) such that
\begin{equation}
\label{eqn T4}
e^{-\lambda_N} \ll N^{-1-\delta}.
\end{equation}

We establish the theorem by showing that, with probability \(1\),
\( \log N \ll r_N(\omega) \ll \log N\) for large \(N\), or equivalently
(in view of Lemmas \ref{Lemma 15} and \ref{Lemma 11})
\begin{equation}
\label{eqn T2}
\lambda'_N \ll r_N(\omega) \ll \lambda'_N
\end{equation}
for \(N > N_0(\omega)\). We apply the Borel-Cantelli
lemma twice to prove that each of the two assertions
of ~(\ref{eqn T2}) holds with probability \(1\). For this
purpose, we must show that if \(C_1, C_2\) are suitably chosen
positive constants, then we have
\begin{equation}
\label{eqn T3}
\sum_{N=0}^{\infty} P(\{ \omega\in\Omega \colon r_N(\omega) > C_1 \lambda'_N \}) < \infty
\end{equation}
and
\begin{equation}
\label{eqn T5}
\sum_{N=0}^{\infty} P(\{ \omega\in\Omega \colon r_N(\omega) < C_2 \lambda'_N \}) < \infty.
\end{equation}
By Lemmas \ref{Lemma 14} and \ref{Lemma 17}, we have
\begin{eqnarray*}
P(\{ \omega\in\Omega\colon r_N(\omega) > C_1 \lambda'_N \}) \leq e^{-\lambda_N} \sum_{d \geq C_1 \lambda'_N} \frac{(\lambda'_N)^d}{d!}
\\
\notag
\leq
e^{-\lambda_N} \left( \frac{e}{C_1} \right)^{C_1 \lambda'_N},
\end{eqnarray*}
provided \(C_1 \geq 1\). Thus, by choosing \(C_1 = e\), we obtain a bound
\( e^{-\lambda_N}\) for the summand of ~(\ref{eqn T3}), and
the inequality ~(\ref{eqn T3}) follows from ~(\ref{eqn T4}).

On the other hand, again by Lemmas \ref{Lemma 14} and \ref{Lemma 17}, we obtain the following estimate
for the summand of ~(\ref{eqn T5}),
\begin{eqnarray*}
P(\{ \omega\in\Omega\colon r_N(\omega) < C_2 \lambda'_N \}) \leq e^{-\lambda_N} \sum_{0 \leq d \leq C_2 \lambda'_N} \frac{(\lambda'_N)^d}{d!}
\\
\notag
\leq
e^{-\lambda_N} \left( \frac{e}{C_2} \right)^{C_2 \lambda'_N},
\end{eqnarray*}
provided \(C_2 \leq 1\). Thus, it suffices to show that
\(C_2\) can be chosen to satisfy, in addition to \(0 < C_2 \leq 1\),
\[
\left( \frac{e}{C_2} \right)^{C_2 \lambda'_N} \ll N^{\delta/2};
\]
for ~(\ref{eqn T5}) will then follow from ~(\ref{eqn T4}).
By Lemmas \ref{Lemma 15} and \ref{Lemma 11}, we know there exists \(D > 0\) such that \(\lambda'_N \leq D \log N \) for
\(N\) sufficiently large. Therefore, we only need to choose
a small positive constant \(C_2\) satisfying
\[
\left( \frac{e}{C_2} \right)^{C_2} \leq e^{\delta/(2 D) },
\]
which is certainly possible since \((e/t)^t \rightarrow 1\) as \(t \rightarrow 0\)
from the positive side. \\

We have now shown that \(\omega\) has each of the desired properties
with probability $1$, and this proves the theorem. \qed

\section{Proof of Theorem \ref{Theorem 0.2}}

Let \(\epsilon > 0\) be given. We define a sequence \(\{ \alpha_j \}\)
by \(\alpha_0 = 1/2\) and
\[
\alpha_j = \frac{1}{2 \ j^{1 - 1/(2 + \epsilon)}}
\]
for \(j \geq 1\). 
It then follows by Lemma \ref{Lemma 11} (with \(\alpha=1/2\), \(c = 1 - 1/(2 + \epsilon)\), and \(c'=0\))
that, with probability \(1\), \(\omega = \{ p_{ b_j }  \}\)
satisfies \(b_j \sim c^* j^{2+ \epsilon} \), where \(c^*\) is some positive constant.

Since the sequence \(\{ \alpha_j \}\) satisfies Hypothesis A,
we have \(\lambda'_N \sim \lambda_N\) by Lemma \ref{Lemma 15}.
Thus, by Lemma \ref{Lemma 11}
we know that there exist positive constants \(D_1\) and \(D_2\) such that
\begin{equation}
\label{eqn Thm 0.2-1}
D_1 q^{- \epsilon n / (2 + \epsilon)} < \lambda_N, \lambda '_N < D_2 q^{- \epsilon n / (2 + \epsilon)}
\end{equation}
for \(N\) sufficiently large.

We again appeal to the Borel-Cantelli Lemma. It follows from this lemma that if
a positive number \(K\) satisfies the property
\begin{equation}
\label{eqn Thm 0.2-2}
\sum_{N=0}^{\infty} P( \{ \omega\in \Omega\colon r_N(\omega) \geq K \}) < \infty,
\end{equation}
then, with probability \(1\), we have
\[
r_N(\omega) < K
\]
for \(N > N_0(\omega)\). \\

We note that, by ~(\ref{eqn Thm 0.2-1}), \(\lambda_N \rightarrow 0\) and
\(\lambda '_N \rightarrow 0\) as \(N \rightarrow \infty\). Thus, by Lemmas \ref{Lemma 14} and \ref{Lemma 17},
we obtain the following estimate for the summand of ~(\ref{eqn Thm 0.2-2}),
\[
P( \{ \omega\in \Omega\colon r_N(\omega) \geq K \}) \leq
e^{-\lambda_N} \sum_{d \geq K} \frac{(\lambda'_N)^d}{d!}
\leq e^{-\lambda_N} \left( \frac{e \lambda'_N}{K} \right)^K
\ll (\lambda'_N)^K
\]
for \(N\) sufficiently large. Since \(q^n \leq N < q^{n+1}\), we have
\[
(\lambda'_N)^K \leq D_2^K q^{- \epsilon n K / (2 + \epsilon)} \ll N^{- \epsilon K/(2 + \epsilon)}.
\]
Therefore, 
provided \(\epsilon K/(2 + \epsilon) > 1\), or equivalently,
\[
K > 1 + 2 \epsilon^{-1},
\]
it is clear that ~(\ref{eqn Thm 0.2-2}) is achieved.
Accordingly we have, with probability \(1\),
\[
r_N(\omega) < 2(1 + \epsilon^{-1})
\]
for \(N > N_1(\epsilon, \omega)\).
This completes the proof of the theorem. \qed

\section{Proof of Theorem \ref{Theorem 0.3}}

Recall we defined \(t_N(\omega)\)
to be
\[
t_N(\omega) = |\{(a,b)\in \mathbb{Z}_{\geq 0} \times \mathbb{Z}_{\geq 0} \colon p_a,p_b \in \omega, \  p_N = p_a - p_b,  \ \deg p_a, \deg p_b \leq \deg p_N \}|.
\]

As before given \(p_N \in \mathbb{F}_q [T] \), we let \(n := n(N) = \deg p_N = \lfloor \log_q N \rfloor\).
It is clear that for $p_N \not = 0$, there exist \(q^{n+1}\) pairs of polynomials \((p_a, p_{b})\)
such that \(p_N = p_a - p_{b}\) and \(\deg p_a\), \(\deg p_{b} \leq n\).
Also, every polynomial of degree less than or equal to \(n\) will appear
as \(p_a\) and \(p_b\) exactly once. Let \(S_{\widehat{u},n}\) denote the set of all
polynomials in $\mathbb{F}_q[T]$ whose degree are less than or equal to \(n\), and the coefficient of $T^n$
is \(\widehat{u} \in \mathbb{F}_q \).
Clearly, we have
\(|S_{\widehat{u},n}| = q^n\). If we consider each polynomial in \(S_{\widehat{u},n}\) as \(p_b\), then the
corresponding set of \(p_a\)'s is \(S_{u,n}\) for some \(u \ne
\widehat{u} \) as \(\deg p_N = n\).

For each \( u \in \mathbb{F}_q \), we consider
\begin{eqnarray*}
t_{N,u }(\omega) = | \{ (a,b)\in \mathbb{Z}_{\geq 0} \times \mathbb{Z}_{\geq 0} \colon p_N = p_a - p_{b}, \text{ where } p_a, p_b \in \omega \text{ and } p_a \in S_{u,n}  \} |.
\end{eqnarray*}

If $p_N = p_a - p_b$, we relabel $p_b$ as $p_{ \widehat{a}}$ to make its correspondence with $p_a$ more explicit.
We form the following two disjoint sets
\[
\mathcal{A}_N = \{a\in \mathbb{Z}_{\geq 0}\colon p_a \in S_{u,n} \} = \{ \iota^{-1} (u) q^n, ..., ( \iota^{-1} (u) +1) q^n-1 \}
\]
and
\[
\widehat{\mathcal{A}}_N = \{\widehat{a}\in \mathbb{Z}_{\geq 0}\colon p_{\widehat{a}} \in S_{\widehat{u},n} \} = \{ \iota^{-1} (\widehat{u}) q^n, ..., ( \iota^{-1} (\widehat{u}) +1) q^n-1 \}.
\]
Let \(M_0:= M_0(N) = |\mathcal{A}_N| = |\widehat{\mathcal{A}}_N| = q^n\).
For convenience, we label the \(M_0\) elements of
\(\mathcal{A}_N\) by \(a_i \  (1 \leq i \leq M_0 ) \),
and the corresponding elements of \( \widehat{\mathcal{A}}_N \) by \(\widehat{a}_i\), in other
words we have $p_N = p_{a_i} - p_{\widehat{a}_i } \ (1 \leq i \leq M_0 )$.

We also define \(\lambda_{N, u }\) and \(\lambda'_{N, u }\) to be
\[
\lambda_{N, u } = \sum_{1 \leq i \leq M_0 } \alpha_{a_i} \alpha_{\widehat{a}_i},
\]
and
\[
\lambda '_{N,  u } = \sum_{1 \leq i \leq M_0 } \frac{ \alpha_{a_i} \alpha_{\widehat{a}_i} }{ 1 - \alpha_{a_i} \alpha_{\widehat{a}_i}}.
\]
With this set up we can recover analogues
of all the previous lemmas in terms of
\(M_0\), \(\lambda_{N, u}\), \(\lambda'_{N, u}\), and \(t_{N,u}(\omega)\),
in place of
\(M\), \(\lambda_{N}\), \(\lambda'_{N}\), and \(r_N(\omega)\), respectively.
Therefore, 
by a similar argument
we obtain Theorem \ref{Theorem 0.2} with \(t_{N,u}(\omega)\) in place of \(r_N(\omega)\).
Since this result holds with probability \(1\), and
\[
t_{N} (\omega) = \sum_{u \in \mathbb{F}_q} t_{N,u}(\omega),
\]
we have our result. \qed

\end{document}